\documentclass[12pt]{amsart}
\usepackage{amssymb}
\numberwithin{equation}{section}

\newtheorem{thm}{Theorem}[section]
\newtheorem{lem}[thm]{Lemma}
\newtheorem{Def}{Definition}[section]
\newtheorem{rem}[Def]{Remark}

\newcommand\calO{{\mathcal{O}}}

\newcommand\calW{{\mathcal{W}}}

\newcommand\bbR{{\mathbb R}}

\newcommand\ol{\overline}

\renewcommand\S{\Sigma}
\newcommand\s{\sigma}
\renewcommand\d{\partial}

\newcommand\D{\nabla}
\newcommand\e{\epsilon}
\renewcommand\b{\beta}

\newcommand\la{\langle}
\newcommand\ra{\rangle}

\newcommand\ric{\operatorname{Ric}}

\newcommand\g{\gamma}
\newcommand\8{\infty}
\renewcommand\a{\alpha}
\newcommand\hess{\operatorname{Hess}}

\renewcommand\th{\theta}
\newcommand\scri{\mathcal{I}}

\newcommand\beq{\begin{eqnarray}}
\newcommand\eeq{\end{eqnarray}}

\begin{document}
\setlength{\baselineskip}{.51cm}
\title[Maximum Principles]
{Maximum Principles for null hypersurfaces and Null Splitting Theorems}

\author{Gregory J. Galloway}
\thanks{Supported in part by NSF grant \#
    DMS-9803566.}
\address{University of Miami, Department of Mathematics and \newline  Computer Science,
Coral Gables, FL, 33124}
\email{galloway@math.miami.edu}


\maketitle

\section{Introduction}

The geometric maximum principle for
smooth (spacelike) hypersurfaces, which is a consequence  
of Alexandrov's \cite{A1} strong maximum
for second order quasi-linear elliptic operators, is a basic tool 
in Riemannian and Lorentzian geometry.  
In \cite{AGH}, extending earlier 
work of Eschenburg \cite{E2}, a version of the 
geometric maximum principle in the Lorentzian setting
was obtained for rough ($C^0$) spacelike hypersurfaces
which obey mean curvature inequalities in the sense
of support hypersurfaces.  In the present paper we 
establish an analogous result for  \emph{null} hypersurfaces (Theorem~\ref{th:3d}) 
and consider some applications.
For the applications, it is important to have a version of the
maximum principle for null hypersurfaces which does not require smoothness.
The reason for this, which is described in more detail in Section 3, is that
the null hypersurfaces which arise most naturally in spacetime geometry
and general relativity, such as black hole event horizons, are 
in general $C^0$ but not $C^1$.  
To establish our basic approach,
we first prove a maximum principle for smooth null hypersurfaces (Theorem~\ref{th:2a}),
and then proceed to the $C^0$ case.  
The general $C^0$ version is then applied to study some rigidity properties of spacetimes
which contain
null lines (inextendible globally maximal null geodesics).
The standard Lorentzian splitting theorem, which is the Lorentzian analogue of the Cheeger-Gromoll
splitting theorem of Riemannian geometry, describes the rigidity of spacetimes which contain
timelike lines (inextendible globally maximal timelike geodesics);
see \cite[Chapter 14]{BEE} for a 
nice exposition. 
Here we show  how the maximum principle for rough null hypersurfaces
can be used to obtain a general ``splitting theorem" for 
spacetimes with null lines (Theorem \ref{th:4a}).  We then consider an 
application  of this null splitting theorem to asymptotically flat 
spacetimes.  We prove
that a nonsingular asymptotically flat (in the sense of Penrose \cite{P1}) vacuum (i.e., Ricci flat)
spacetime which contains a null line is isometric to Minkowski space (Theorem~\ref{th:4c}).

In Section 2 we review the relevant aspects of the geometry of null hypersurfaces
and present the maximum principle for smooth 
null hypersurfaces.  In Section 3  we present the maximum 
principle for $C^0$ null hypersurfaces. In Section 4 we consider the
aforementioned applications.  For basic notions used below from
Lorentzian geometry and causal theory, 
we refer the reader to the excellent references,
\cite{BEE}, \cite{HE}, \cite{ON} and \cite{P2}.

\section{The maximum principle for smooth null hypersurfaces}
\subsection{The geometry of null hypersurfaces}
Here we review some aspects of the geometry of null hypersurfaces.
For further details, see e.g. \cite{K} which is written from a similar
point of view.

Let $M$ be a spacetime, i.e., a smooth time-oriented Lorentzian manifold, of dimension $n\ge 3$.  
We denote the Lorentzian metric on $M$ by $g$ or $\la\,,\ra$.  
A {\it (smooth) null hypersurface} in $M$
is a smooth co-dimension one embedded submanifold $S$ of $M$ such that the pullback of the metric
$g$ to $S$ is degenerate.  
Because of the Lorentz signature of $g$, 
the null space of the pullback is one dimensional at each point of $S$. 
Hence, every null hypersurface $S$ admits a smooth nonvanishing future directed
null vector field $K\in \Gamma TS$ such that the normal space of $K$ at $p\in S$ coincides
with the tangent space of $S$ at $p$, i.e., $K_p^{\perp} = T_pS$ for all $p\in S$.  
It follows, in
particular, that tangent vectors to $S$ not parallel to $K$ are spacelike.  
It is well-known that the integral curves of $K$, when suitably parameterized, are null geodesics.
These integral curves are called the {\it null geodesic generators\/} of $S$.
We note that the vector field $K$ is unique up to a positive (pointwise) scale factor. 

Since $K$ is orthogonal to $S$ we can introduce the null Weingarten map and null second fundamental
form of $S$ with respect $K$ in a manner roughly analogous to what is done for spacelike
hypersurfaces  or hypersurfaces in a Riemannian manifold.  

We introduce the following equivalence
relation on tangent vectors:  For  $X, X'\in T_pS$, $X'=X \mbox{ mod } K$ if and only if 
$X' - X = \lambda K$ for some $\lambda \in \Bbb R$.  Let $\ol X$ denote the equivalence class
of $X$.  Simple computations show that if
$X'=X \mbox{ mod } K$ and $Y'=Y \mbox{ mod } K$  then $\la X',Y'\ra = \la X,Y\ra$ and 
$\la \D_{X'} K,Y'\ra = \la \D_X K,Y\ra$, where $\D$ is the Levi-Civita connection of $M$.
Hence, for various quantities of interest, components along $K$ are not of interest.  For this
reason one works with the tangent space of $S$ modded out by $K$, i.e., let $T_pS/K = \{\ol X: X\in
T_pS\}$ and $TS/K = \cup_{p\in S}T_pS/K$.  $TS/K$ is a rank $n-2$ vector bundle over $S$.  This
vector bundle does not depend on the particular choice of null vector field $K$.  
There is a natural positive definite metric $h$ in $TS/K$ induced from $\la\,,\ra$: For each $p\in S$, 
define $h:T_pS/K\times T_pS/K\to \Bbb R$  by $h(\ol X,\ol Y) = \la X,Y\ra$. From remarks
above, $h$ is well-defined.  

The {\it null Weingarten map\/} $b=b_K$ of $S$ with respect to $K$ is, for each point $p\in S$, a
linear map $b: T_pS/K\to T_pS/K$ defined by $b(\ol X) = \ol{\D_X K}$.  It is easily verified
that $b$ is well-defined.  Note if $\widetilde K = fK$, $f\in C^{\8}(S)$, is any other future directed
null vector field tangent to $S$, then $\D_X \widetilde K= f\D_X K \mbox{ mod } K$.  It follows that
the Weingarten map $b$ of $S$ is unique up to positive scale factor and that $b$ at a given
point $p\in S$ depends only on the value of $K$ at $p$.    

A standard computation shows, $h(b(\ol X), \ol Y) = \la \D_X K, Y\ra =
\la X, \D_Y K\ra = h(\ol X,b(\ol Y))$.  Hence $b$ is self-adjoint with respect to $h$.
The {\it null second fundamental form\/} $B=B_K$ of $S$ with respect to $K$ is the bilinear form
associated to $b$ via $h$:  For each $p\in S$, $B:T_pS/K\times T_pS/K\to \Bbb R$ is defined by
$B(\ol X,\ol Y) = h(b(\ol X),\ol Y)= \la\D_X K,Y\ra$.  Since $b$ is self-adjoint, $B$ is
symmetric. 
We say that $S$ is \emph{totally geodesic\/} iff $B\equiv 0$.  This has the usual
geometric consequence: A geodesic in $M$ starting tangent to a totally 
geodesic null hypersurface $S$ remains in $S$.  Null hyperplanes in Minkowski space
are totally geodesic, as is the event horizon in Shwarzschild spacetime.     

The {\it null mean curvature\/}  of $S$ with respect to $K$
is the smooth scalar field $\theta\in C^{\8}(S)$ defined by $\theta ={\rm tr}\,b$.  
Let $e_1,e_2, ... e_{n-2}$ be $n-2$ orthonormal spacelike vectors
(with respect to $\la,\ra$) tangent to $S$ at
$p$. Then $\{\ol e_1,\ol e_2, ... \ol e_{n-2}\}$ is an orthonormal basis (with respect to $h$) of
$T_pS/K$.   Hence at $p$, 
\beq
\theta 
& = & {\rm tr}\,b = \sum_{i=1}^{n-2} h(b(\ol e_i),\ol e_i) \nonumber\\
& = &  \sum_{i=1}^{n-2} \la \D_{e_i}K,e_i\ra.  \label{eq:2a}
\eeq

Let $\S$ be the {\it properly transverse\/} intersection of a hypersurface $P$ in $M$ with $S$.
By {\it properly transverse\/} we mean that $K$ is not tangent to $P$ at any point of $\S$. 
Then $\S$ is a smooth $(n-2)$-dimensional spacelike submanifold of $M$ contained
in $S$ which meets $K$ orthogonally. From Equation \ref{eq:2a},
$\theta|_{\S} = {\rm div}_{\S}K$, and hence the null mean curvature gives a measure of the
divergence of the null generators of $S$. Note that
if
$\widetilde K = fK$ then
$\widetilde\theta =f\theta$.  Thus the null mean curvature inequalities $\theta\ge 0$, $\theta\le 0$,
are invariant under positive rescaling of $K$.  In Minkowski space, a future null cone 
$S =\d I^+(p)\setminus\{p\}$  (resp., past null cone $S =\d I^-(p) \setminus\{p\}$) has positive null mean
curvature,
$\theta >0$  (resp.,  negative null mean curvature, $\theta <0$).  

The  null second fundamental form of a null hypersurface obeys a well-defined 
comparison theory roughly similar to the comparison theory satisfied by the second
fundamental forms of a family of parallel spacelike hypersurfaces (cf., Eschenburg \cite{E1},
which we follow in spirit).  

Let $\eta :(a,b)\to M$, $s\to \eta(s)$, be a future directed affinely parameterized 
null geodesic generator of $S$.  For each $s\in (a,b)$, let 
$$
b(s) = b_{\eta'(s)}: T_{\eta(s)}S/\eta'(s) \to T_{\eta(s)}S/\eta'(s)
$$
be the Weingarten map based at $\eta(s)$ with respect to the null vector $K=\eta'(s)$. 
The one parameter family of Weingarten maps $s\to b(s)$,
obeys the following Ricatti equation,
\begin{eqnarray}\label{eq:2b}
b'+b^2 +R = 0 .  
\end{eqnarray} 
Here $'$ denotes covariant differentiation in the direction $\eta'(s)$; if $X=X(s)$ is a
vector field along $\eta$ tangent to $S$, we define,
\beq\label{eq:2bb}
b'(\ol X) = b(\ol X)'- b(\ol{X'}).  
\eeq
$R:T_{\eta(s)}S/\eta'(s)\to T_{\eta(s)}S/\eta'(s)$ is the curvature endomorphism defined by
$R(\ol X) = \ol{R(X,\eta'(s))\eta'(s)}$, where $(X,Y,Z)\to R(X,Y)Z$ is the Riemann
curvature tensor of $M$, $R(X,Y)Z = \D_X\D_YZ-\D_Y\D_XZ - \D_{[X,Y]}Z$.

We indicate the proof of Equation \ref{eq:2b}.  Fix a point $p = \eta(s_0)$,
$s_0\in (a,b)$, on $\eta$.  On a neighborhood $U$ of $p$ in $S$ we can 
scale the null vector field $K$ so that 
$K$ is a geodesic 
vector field, $\nabla_KK=0$, and so that $K$, restricted to $\eta$, is
the velocity vector field to $\eta$, i.e., for each $s$ near $s_0$, $K_{\eta(s)}=
\eta'(s)$. Let $X\in T_pM$.  Shrinking $U$ if necessary, we can extend $X$ to a
smooth vector field on $U$ so that $[X,K] =\nabla_XK - \nabla_KX = 0$. 
Then, $R(X,K)K = \nabla_X\nabla_KK-\nabla_K\nabla_XK -\nabla_{[X,K]}K =
-\nabla_K\nabla_KX$.  Hence along $\eta$ we have, $X'' = -R(X,\eta')\eta'$
(which implies that $X$, restricted to $\eta$, is a Jacobi field along
$\eta$).  Thus, from Equation \ref{eq:2bb}, at the point $p$ we have,
\begin{eqnarray*}
b'(\overline X) & = & \overline{\nabla_XK}\,' - b(\overline{\nabla_KX})   = 
\overline{\nabla_KX}\,'-b(\overline{\nabla_XK})\\ 
   & = &  \overline{X''} - b(b(\overline X))  =  - \overline{R(X,\eta')\eta'} - b^2(\overline X) \\
& = & - R(\overline X) - b^2(\overline X),
\end{eqnarray*}
which establishes Equation \ref{eq:2b}.

By taking the trace of \ref{eq:2b} we obtain the following formula for the derivative
of the null mean curvature $\theta=\theta(s)$ along $\eta$,
\begin{eqnarray}
\theta' = -{\rm Ric}(\eta',\eta') - \sigma^2 - \frac1{n-2}\theta^2,  \label{eq:2c}
\end{eqnarray}
where $\sigma$, the shear scalar, is the trace of the square of the trace free part of
$b$. Equation \ref{eq:2c} is the well-known  Raychaudhuri equation (for an
irrotational null geodesic congruence) of relativity.   This equation shows how
the Ricci curvature of spacetime influences the null mean curvature of a null hypersurface.

\subsection{The maximum principle for smooth null hypersurfaces}
The aim here is to prove the geometric maximum principle for smooth
null hypersurfaces.  Because of its natural invariance, we restrict attention
to the zero mean curvature case.
In the statement we make use of the following intuitive terminology.  Let $S_1$ and
$S_2$  be  null hypersurfaces that meet at a point $p$.  We say that $S_2$ \emph{lies to the future 
(resp., past) side of $S_1$ near $p$} provided for some neighborhood $U$ of $p$ in $M$ in which $S_1$ is
closed and achronal,
$S_2\cap U \subset J^+(S_1\cap U,U)$ (resp., $S_2\cap U \subset J^-(S_1\cap U,U)$).

\begin{thm}\label{th:2a} Let $S_1$ and $S_2$ be smooth null hypersurfaces in a spacetime $M$.
Suppose,
\begin{enumerate}
\item[(1)] $S_1$ and $S_2$ meet at $p\in M$ and $S_2$ lies to the future side of $S_1$ near $p$,
and
\item[(2)] the null mean curvature scalars $\theta_1$ of $S_1$, and $\theta_2$ of $S_2$,
satisfy, \hfil
\linebreak$\theta_1 \le 0 \le \theta_2$.  
\end{enumerate}
Then $S_1$ and $S_2$ coincide near $p$ and this common null hypersurface has 
null mean curvature $\theta = 0$. 
\end{thm}

The proof is an application of Alexandrov's \cite{A1} strong maximum principle 
for second order quasi-linear elliptic operators.  It will be convenient to state the precise
form of this result needed here.  

For connected open sets $\Omega\subset \bbR^n$
and $U\subset \bbR^n\times \bbR \times \bbR^n$, we say $u\in C^2(\Omega)$ is $U$-{\it
admissible\/} 
provided $(x,u(x),\d u(x))\in U$
for all $x=(x^1, x^2,..., x^n)\in \Omega$, where  $\d u = (\d_1u,\d_2u,...,\d_nu)$, and 
$\d_iu =\frac{\d u}{\d x^i}$.  

Let $Q=Q(u)$ be a second order quasi-linear operator, i.e., for $U$-admissible $u\in C^2(\Omega)$, consider
\begin{eqnarray}\label{eq:2d}
Q(u) = \sum_{i,j=1}^n a^{ij}(x,u,\d u)\d_{ij}u + b(x,u,\d u), 
\end{eqnarray}
where $a^{ij}, b\in C^1(U)$, $a^{ij} = a^{ji}$, and $\d_{ij} = \frac{\d^2}{\d u^j\d u^i}$.
The operator $Q$ is {\it elliptic\/} provided for each $(x,r,p)\in U$, and for all
$\xi = (\xi^1, ...\xi^n)\in \bbR^n$, $\xi \ne 0$,
$$
\sum_{i,j=1}^n a^{ij}(x,r,p)\xi^i\xi^j > 0.
$$
We now state the form of the strong maximum principle for second order quasi-linear
elliptic operators most suitable for our purposes.

\begin{thm}\label{th:2b} 
Let $Q=Q(u)$ be a second order quasi-linear elliptic operator as described above.
Suppose the $U$-admissible functions $u,v\in C^2(\Omega)$ satisfy, 
\begin{enumerate}
\item[(1)] $u\le v$ on $\Omega$ and $u(x_0)=v(x_0)$ for some $x_0\in \Omega$, and 
\item[(2)] $Q(v) \le Q(u)$ on $\Omega$.
\end{enumerate}
Then $u=v$ on $\Omega$.
\end{thm}

The idea of the proof of Theorem \ref{th:2a} is to intersect, in a properly transverse
manner, the null hypersurfaces
$S_1$ and
$S_2$ with a {\it timelike\/} (i.e., Lorentzian in the induced metric)  hypersurface through the point 
$p$, and to show that the  spacelike intersections agree.  Analytically, intersecting the null hypersurfaces in
this manner  reduces the problem to a nondegenerate elliptic one.  In order to apply Theorem \ref{th:2b}
we need a suitable analytic expression  for the null mean curvature, which we now derive.

Let $p$ be a point in a spacetime $M$, and let $P$ be a timelike hypersurface passing through 
$p$. Let $V$ be a connected spacelike hypersurface
in $P$ (and hence a co-dimension two spacelike submanifold of $M$) passing through $p$.
Via the normal exponential map along $V$ in $P$, we can assume, by shrinking $P$ if necessary,
that $P$ can be expressed as,
\begin{eqnarray}\label{eq:2e}
P= (-a,a)\times V, 
\end{eqnarray}
and that the induce metric on $P$ takes the form,
\begin{eqnarray}\label{eq:2f}
ds^2= -dt^2 + \sum_{i=1}^{n-2}g_{ij}(t,x)dx^idx^j,
\end{eqnarray}
where $x=(x_1,...,x_{n-2})$ are coordinates in $V$ centered on $p$. 

Let $S$ be a null hypersurface which meets $P$ properly transversely in a spacelike
hypersurface $\S$ in $P$.  By adjusting the size of $P$ and $S$ if necessary, we may
assume that $\S$ can be expressed as a graph over $V$.  Thus, there exists $u\in C^{\8}(V)$
such that 
$\S={\rm graph}\,u = \{(u(x),x)\in P: x\in V\}$.
Let $H(u)$ denote the mean curvature of $\S= {\rm graph}\,u$.  (By our sign conventions
the mean curvature of $\S$ is $+$ the divergence of the future pointing normal along $\S$.)  To describe 
$H(u)$ we introduce the following notation.  Let $h$ be the Riemannian metric on $V$ whose
components are given by $h_{ij}(x)= g_{ij}(u(x),x)$, and let $h^{ij}$ be the $i,j$th entry of the
inverse  matrix $[h_{ij}]^{-1}$.  Let $\D u$ denote the gradient of $u$.  In terms of coordinates,
$\D u = \sum_j u^j\d_j$, where $u^j = \sum_i h^{ij}\d_i u$.  Finally, introduce the quantity,
\beq
\nu  =  \frac{1}{\sqrt{1-|\D u|^2}}.   \nonumber
\eeq
The positivity of $\nu$ is equivalent to  $\S={\rm graph}\,u$ being spacelike.  With respect to
these  quantities, we have (cf.,  \cite{AGH}),
\beq
H(u) = \sum_{i,j=1}^{n-2} a^{ij}(x,u,\d u)\d_{ij} u + b(x,u,\d u), \label{eq:2g}
\eeq
where $a^{ij}= \nu h^{ij} + \nu^3u^iu^j$ and $b$ is a polynomial expression in $\d_iu$, $h_{ij}$,
$h^{ij}$, $\d_t g_{ij}(u(x),x)$ 
and $\nu$.  From the form of $a^{ij}$, it is clear that $H = H(u)$ is a second
order quasi-linear elliptic operator.  

Let $K$ be a future directed null vector field on $S$.  Since $K$ is orthogonal to $\S$, by rescaling
we may assume that along $\S$ , $K=Z+N$, where 
$Z$ is the future directed unit normal vector field to $\S$ in $P$ and
$N$ is one of the two unit spacelike normal vector fields to $P$ in $M$.
Let $\theta$ be the null mean curvature of $S$ with respect to $K$.  We obtain an expression for
$\theta$ along $\S$.  Let $B_P$ denote the second fundamental form of $P$ with respect to $N$,
and let $B_{\S}$ denote the second fundamental form of $\S$ in $P$ with respect to $Z$.  Then
for $q \in \S$ and vectors $X,Y\in T_q\S$, $B_P(X,Y) = \la\D_X N,Y\ra$, and $B_{\S}(X,Y)=
\la\overline{\D}_X Z,Y\ra = \la\D_X Z,Y\ra$, where $\overline{\D}$ is the induced Levi-Civita
connection on $P$.

Now let $\{e_1, e_2, ..., e_{n-2}\}$ be an orthonormal basis for $T_q\S$.  Then the value of
$\theta$ at $q$ is given by,
\beq
\theta & = & \sum_{i=1}^{n-2} \la \D_{e_i}K,e_i\ra = 
\sum_{i=1}^{n-2} \la \D_{e_i}Z,e_i\ra + \sum_{i=1}^{n-2} \la \D_{e_i}N,e_i \ra \nonumber \\
& = & \sum_{i=1}^{n-2} B_{\S}(e_i,e_i) + \sum_{i=1}^{n-2} B_{P}(e_i,e_i) \nonumber \\
& = & H_{\S} + B_P(Z,Z) + H_P, \label{eq:2h}
\eeq
where $H_{\S}$ is the mean curvature of $\S$ and $H_P$ is the mean curvature of $P$.  
In the notation introduced above, 
\beq
Z & = & \nu(\d_0 + \D u) \label{eq:2hh} \\
& = & \sum_{i=0}^{n-2}\nu u^i\d_i,  \nonumber
\eeq 
where $\d_0= \frac{\d}{\d t}$, $u^0 = 1$, and as above, $u^i= \sum_{j=1}^{n-2}h^{ij}\d_j u$, 
$i= 1, ..., n-2$. Hence,
\beq
B_P(Z,Z) & = & B_P(\sum_{i=0}^{n-2}\nu u^i\d_i,\sum_{i=0}^{n-2}\nu u^i\d_i)  \nonumber\\
& = &  \sum_{i,j=1}^{n-2} \nu^2\beta_{ij}(u)u^iu^j, \label{eq:2i}
\eeq
where for $x\in V$, $\beta_{ij}(u)(x) = B_P(\d_i,\d_j)|_{(u(x),x)}$.

Now let $\theta(u)$ denote the null mean curvature of $S$ along $\S={\rm graph}\,u$.  Equations
\ref{eq:2h} and \ref{eq:2i} give,
$$
\theta(u) = H(u) + \sum_{i,j=1}^{n-2} \nu^2\beta_{ij}(u)u^iu^j + \alpha(u),
$$
where $\a(u)$ is the function on $V$ defined by $\a(u)(x) = H_P(u(x),x)$.  Thus, by \ref{eq:2g},
we finally arrive at,
\beq
\theta(u) = \sum_{i,j=1}^{n-2} a^{ij}(x,u,\d u)\d_{ij} u + b_1(x,u,\d u) \label{eq:2j}
\eeq
where,
\beq
b_1(x,u,\d u) = b(x,u,\d u) + \sum_{i,j=1}^{n-2} \nu^2\beta_{ij}(u)u^iu^j + \alpha(u),\label{eq:2jj}
\eeq
and where $a^{ij}$ and $b$ are as in Equation \ref{eq:2g}.  In particular, $\theta = \theta(u)$
is a second order quasi-linear elliptic operator with the same symbol as the mean curvature
operator for spacelike hypersurfaces in $P$.

\smallskip\noindent
{\it Proof of Theorem 2.1:}  Let $P$ be a timelike hypersurface passing through $p$, as described
by equations \ref{eq:2e} and \ref{eq:2f}.  $S_1$ and $S_2$ have a common null tangent at $p$.
Choose $P$ so that it is transverse to this tangent.  Then, by choosing $P$ small enough 
the intersections $\S_1 = S_1\cap P$ and $\S_2 = S_2\cap P$ will be properly transverse,
and hence $\S_1$ and $\S_2$ will be spacelike hypersurfaces
in $P$.  Let $K_i$, $i=1,2$, be the null vector field on $S_i$ with respect to which the
null mean curvature function $\theta_i$ is defined.  Let $N$ be the unit normal to $P$
pointing to the same side of $P$ as $K_1|_{\S_1}$ and $K_2|_{\S_2}$.  By rescaling we can
assume $K_i|_{\S_i} = Z_i+ N|_{\S_i}$, where $Z_i$ is the future directed unit normal to $\S_i$ in $P$.

Let $u_i=u_i(x)$, $i=1,2$, be the smooth function on $V$ such that $\S_i= {\rm graph}\,u_i$.
From the hypotheses of Theorem \ref{th:2a} we know, 
\begin{enumerate}
\item[(1)] $u_1\le u_2$ on $V$ and $u_1(p)=u_2(p)$, and 
\item[(2)] $\theta(u_2) \le\theta(u_1)$ on $V$.
\end{enumerate} 

By Theorem \ref{th:2b}, we conclude that $u_1=u_2$ on $V$, i.e. $\S_1=\S_2$. Now,
$S_i$, $i=1,2$, is obtained locally by exponentiating out from $\S_i$ in the 
orthogonal direction $K_i|_{\S_i} = Z_i+ N|_{\S_i}$.  It follows that $S_1$ and $S_2$ agree near $p$, i.e.,
there is a spacetime neighborhood $\calO$ of $p$ such that $S_1\cap\calO =S_2\cap\calO = S$,
and $S$ has null mean curvature $\theta = 0$.        
\hfill\endproof
\smallskip

\section{The maximum principle for $C^0$ null hypersurfaces}
\subsection{$C^0$ null hypersurfaces}

The usefulness of the maximum principle for smooth null hypersurfaces obtained in 
the previous section
is  limited by the fact that the most interesting null hypersurfaces arising
in general relativity, such as black
hole event horizons and Cauchy horizons, are rough, i.e., are $C^0$ but in general not $C^1$.
The aim of this section is to present a maximum principle for rough null hypersurfaces,
similar in spirit to the maximum principle for rough spacelike hypersurfaces obtained
in \cite{AGH}.   

Horizons and other null hypersurfaces commonly occurring in relativity arise essentially as the null
portions of
\emph{achronal boundaries} which are sets of the form $\d I^{\pm}(A)$, $A\subset M$.
Achronal boundaries are always $C^0$ hypersurfaces, but simple examples illustrate that they 
(and their null portions) may fail
to be differentiable at certain points.  Consider, for example, the set 
$S = \d I^{-}(A) \setminus~A$, where $A$ consists of two disjoint closed disks in the $t=0$
slice of Minkowski $3$-space.
This surface, which represents the merger of two
truncated cones, has a ``crease", i.e., a curve of nondifferentiable points (corresponding to the 
intersection of the two cones) but which otherwise is
a smooth null hypersurface. 

An important feature of the null portion of an achronal boundary is that it is
ruled by null geodesics 
which are either past or future inextendible within the hypersurface.
This is illustrated in the example above.  $S$ is ruled by null geodesics
which are future inextendible in $S$, but which are in general not past inextendible
in $S$.  Null geodesics in $S$ that meet the crease when extended toward the past
leave $S$  when extended further, and hence have past end points on $S$.

We now formulate a general definition of $C^0$ null hypersurface which captures the essential
features of these examples. 

A set $A\subset M$ is said to be achronal if no two points
of $A$ can be joined by a timelike curve.  
$A\subset M$ is locally achronal if for each $p\in A$
there is a neighborhood $U$ of $p$ such that $A\cap U $ is achronal in $U$.         
A \emph{$C^0$ nontimelike hypersurface} in $M$ is a topological hypersurface $S$ in $M$ which is
locally achronal.  We remark that for each  $p \in S$, there exist arbitrarily small
neighborhoods $U$ of $p$ such that
$S\cap U$ is closed and achronal in $U$, and  for each $q\in U\setminus S$, 
either $q\in I^+(S\cap U,U)$ or $q\in I^-(S\cap U,U)$.

\begin{Def}\label{def:3a} A $C^0$ \emph{future null hypersurface\/} in $M$ is a
nontimelike hypersurface $S$ in $M$ such that
for each $p\in S$ and any neighborhood  $U$ of $p$ in 
which $S$ is achronal, there exists a point $q\in S$, $q\ne p$, such that $q \in
J^+(p,U)$. 
\end{Def}

Since $q\in J^+(p,U)\setminus I^+(p,U)$, there is a null geodesic segment $\eta$ from $p$ to 
$q$. The segment $\eta$ must be contained in $S$,  for otherwise at some point $\eta$ would enter
either $I^+(S,U)$ or $I^-(S,U)$, which would contradict the achronality of $S$ in $U$.
The geodesic $\eta$  can be extended further to the future in $S$: Choose $r$ in
$S \cap J^+(q,U)$, $r \ne q$. The null geodesic from $q$ to $r$ in $S$ must smoothly extend the one 
from $p$ to $q$, otherwise there would be an achronality violation of $S$ in $U$.  Thus, for each
point $p \in S$ there is a future directed null geodesic in $S$ starting at, or passing through, $p$ 
which is future inextendible in $S$, i.e., which does not have a future end point in $S$.  These null
geodesics are called the \emph{null geodesic generators} of
$S$.  They may or may not have past end points in $S$.  Summarizing, a $C^0$ future null
hypersurface is a locally achronal topological hypersurface $S$ of $M$ which is ruled by future
inextendible null geodesics.  A $C^0$ past null hypersurface is defined in a time-dual
manner. 

Let $S$ be a $C^0$ future null hypersurface.  Adopting 
the terminology introduced in \cite{CG} for event horizons,  a \emph{semi-tangent} of $S$ 
is a future directed null vector $K$ which is tangent to a null generator
of $S$.  We do not want to distinguish between semi-tangents based at the same point and
pointing in the same null direction, so we assume the  semi-tangents of $S$ have been 
uniformly normalized in some manner, e.g., by requiring each semi-tangent to have unit
length with respect to some auxilliary Riemannian metric on $M$.  Then note that the local
achronality of $S$ implies that at each interior point (non-past end point) of a 
null generator of $S$ there is a unique semi-tangent at that point.  Techniques from
\cite{CG} can be adapted to prove the following.  

\begin{lem}\label{th:3a} Let $S$ be a $C^0$ future null hypersurface in a spacetime $M$.
\begin{enumerate}
\item If $p_n \rightarrow p$ in $S$ and $X_n \rightarrow X$ in $TM$, where,
for each $n$, $X_n$ is a semi-tangent of $S$ at $p_n$ then $X$ is a semi-tangent of $S$.
\item Suppose $p$ is an interior point of a null generator of $S$, and let $X$ be the
unique semi-tangent of $S$ at $p$.  Then semi-tangents of $S$ at points near $p$
must be close to $X$, i.e., if $X_n$ is any 
semi-tangent of $S$ at $p_n$, and $p_n \rightarrow p$ then $X_n \rightarrow X$. 
\end{enumerate}
\end{lem}

The proof of the maximum principle for $C^0$ null hypersurfaces will proceed in a fashion
similar to the smooth case.  Thus we will need to consider the intersection of a 
$C^0$ null hypersurface $S$ with a smooth timelike hypersurface $P$.  

\begin{lem}\label{th:3b} Let $S$ be a $C^0$ future null hypersurface and let $p\in S$ be an
interior point of a null generator $\eta$ of $S$.  Let $P$ be a smooth
timelike hypersurface passing through $p$ transverse to $\eta$.  Then there exists a
neighborhood $\calO$ of $p$ in $P$ such that $\S = S\cap P$ is a partial Cauchy surface
in $\calO$, i.e., $\S$ is a closed acausal $C^0$ hypersurface in $\calO$.    
\end{lem} 

\begin{proof}  The proof uses the edge concept, in particular the result that an achronal
set is a closed $C^0$ hypersurface if and only if it is edgeless; see e.g., 
Corollary~26, p. 414 in \cite{ON}. 
Let $U$ be a neighborhood of $p$ in $M$ in which $S$ is achronal and edgeless.
Then $V =U\cap P$ is a neighborhood of $p$ in $P$ in which $\S$ is achronal and edgeless in
$P$.  Hence, $\S$ is a closed achronal $C^0$ hypersurface in $V$, and it remains to show that it is
actually acausal in a perhaps smaller neighborhood.
Suppose there exists a sequence of neighborhoods $V_n\subset V$ of p,
which shrink to $p$, such that $\S$ is not acausal in $V_n$ for each $n$.  Then, for each $n$,
there exists a pair of points $p_n,q_n\in \S\cap V_n$ such that $p_n \to p$ and $q_n \in J^+(p_n, V_n)$.
Hence for each $n$, there exists a $P$-null geodesic $\eta_n$ from $p_n$ to $q_n$. 
Now, $\eta_n$ is a causal curve in $U$, and, in fact,  must be a null geodesic in $U$,
since otherwise we would have $q_n \in I^+(p_n,U)$, which would violate the achronality of
$S$ in $U$. Hence $\eta_n\subset S$, and the initial tangent $X_n$ to $\eta_n$, when suitably normalized,
is   a semi-tangent of $S$ at $p_n$.  By the second part of Lemma \ref{th:3a}, $X_n \to X$, where $X$
is tangent to $\eta$ at $p$.  But $X$ is tangent to $P$, since each $X_n$ is, which 
contradicts the assumption that $P$ is transverse to $\eta$ at $p$.   
\end{proof}

We now extend the meaning of mean curvature inequalities to $C^0$ null hypersurfaces. 
The idea, motivated by previous work involving spacelike hypersurfaces (\cite{E2}, \cite{AGH})
is to use smooth null support hypersurfaces.  \emph{Henceforth we set the scale for all null vectors
on $M$ by requiring that they have unit length with respect to a fixed Riemannian metric on $M$}.

\begin{Def}\label{def:3b} Let $S$ be a $C^0$ future null hypersurface in $M$.  We say that $S$ has null
mean curvature $\th \ge 0$ \textup{in the sense of support hypersurfaces} provided 
for each $p\in S$ and for each $\e >0$ there exists a smooth (at
least $C^2$)  null hypersurface $S_{p,\e}$ such that,
\begin{enumerate}
\item[(1)] $S_{p,\e}$ is a past support hypersurface for $S$ at $p$, i.e., $S_{p,\e}$ passes 
through $p$ and lies to the past side of $S$ near $p$, and
\item[(2)] the null mean curvature of $S_{p,\e}$ at $p$ satisfies $\th_{p,\e}\ge -\e$. 
\end{enumerate}
\end{Def}

For example, if $p$ is a point in Minkowski space, the future null cone $\d I^+(p)$ has null
mean curvature $\th \ge 0$ in the sense of support hypersurfaces.  One can use null hyperplanes,
even at the vertex, as support hypersurfaces.  Another, less trivial example is that of a
black hole event horizon $H = \d I^-(\scri^+)$ in an asymptotically flat black hole 
spacetime $M$.  Here $\scri^+$ refers to future null infinity; see Section~4.  Assuming the
generators of $H$ are future complete and  $M$ obeys the null energy condition, $\ric(X,X)\ge 0$,
for all null vectors $X$, it follows from Lemma \ref{th:4b} in Section 4 that $H$ has null
mean curvature $\theta \ge 0$ in the sense of support hypersurfaces.  
This observation and other consequences of the existence of smooth null
support hypersurfaces provided the initial impetus for the development of a proof
of the black hole area theorem under natural regularity conditions, i.e. the
regularity implicit in the fact that $H$ is an achronal boundary, cf. \cite{C+}.     

With the notation as in Definition \ref{def:3b}, let $B_{p,\e}$ denote the null second
fundamental form of $S_{p,\e}$ at $p$.  We say that the collection of null second fundamental
forms $\{B_{p,\e}:p\in S,\e>0\}$ is \emph{locally bounded from below} provided that for all
$p\in S$ there is a neighborhood $\calW$ of $p$ in $S$ and a constant $k>0$ such that
\beq\label{eq:3bb}
B_{q,\e} \ge -k h_{q,\e}\quad\mbox{for all } q\in \calW \mbox{ and } \e >0,
\eeq     
where $h_{q,\e}$ is the Riemannian metric on $T_qS_{q,\e}/K_{q,\e}$, as defined in Section 2.1.  
This technical condition arises in the statement of the maximum principle for $C^0$
null hypersurfaces, and is satisfied in many natural geometric situations for essentially a priori
reasons.  

\begin{lem}\label{th:3c} Let $S$ be a $C^0$ future null hypersurface and let $W$ be
a smooth null hypersurface which is a past support hypersurface for $S$ at $p$.  If
$K\in T_pW$ is the future directed (normalized) null tangent of $W$ at $p$ then
$K$ is a semi-tangent of $S$ at $p$.     
\end{lem}

\begin{proof}    
Let $U$ be a neighborhood of $p$ such that $S\cap U$ is closed and achronal 
in $U$ and $W\cap U \subset J^-(S\cap U,U)$.  
For simplicity, we may assume that $S\subset U$ and $W \subset J^-(S,U)$. 
Let $\eta\subset U$ be an initial segment in $U$ of the future directed null generator of $W$
starting at $p$ with initial tangent $K$.  Then $\eta \subset J^-(S,U)\cap J^+(S,U)$. 
By the remark before Definition \ref{def:3a}, if $\eta$ leaves $S$ at some point, it will
enter either $I^-(S,U)$ or $I^+(S,U)$.  Either case leads to an achronality violation. Hence,
$\eta$ must be a null generator of $S$, which implies that $K$ is a semi-tangent of $S$.
\end{proof}

If $S$ is a $C^0$ past null hypersurface, one defines in a time-dual fashion what it means for 
$S$ to have null mean curvature $\theta\le 0$ in the sense of support hypersurfaces.  In this case
one uses smooth null hypersurfaces which lie locally to the \emph{future} of $S$.  Although, in
principle, one can also consider future null hypersurfaces with nonpositive null mean curvature, this
appears to be a less useful notion, as future support hypersurfaces cannot 
typically be constructed at past end points of generators.

\subsection{The maximum principle for $C^0$ null hypersurfaces}

The aim now is to present a proof of the geometric maximum principle stated below.  Unless otherwise
stated, we continue to assume that all null vectors are normalized to unit length with respect to
a fixed background Riemannian metric.  

\begin{thm}\label{th:3d} Let $S_1$ be a $C^0$ future null hypersurface and let $S_2$ be a $C^0$ past null
hypersurface in a spacetime $M$. Suppose,
\begin{enumerate}
\item[(1)] $S_1$ and $S_2$ meet at $p\in M$ and $S_2$ lies to the future side of $S_1$ near $p$,
\item[(2)] $S_1$ has null mean curvature $\th_1 \ge 0$ in the sense of support hypersurfaces,
with null second fundamental forms $\{B_{p,\e}:p\in S_1,\e>0\}$ locally bounded from below, and 
\item[(3)] $S_2$ has null mean curvature $\th_2\le 0$ in the sense of support hypersurfaces.
\end{enumerate}
Then $S_1$ and $S_2$ coincide near $p$, i.e., there is a neighborhood $\calO$of $p$ such that
$S_1\cap \calO =  S_2\cap \calO$.  Moreover, $S_1\cap \calO =  S_2\cap \calO$ is a smooth
null hypersurface with null mean curvature $\th = 0$.   
\end{thm}

The proof proceeds in a similar fashion to the proof of Theorem \ref{th:2a}.  Instead of Theorem 
\ref{th:2b}, we use the strong maximum principle for weak (in the sense of support functions)
sub and super solutions of second order quasi-linear elliptic equations obtained in
\cite{AGH}.  We will state here a restricted form of this result, adapted to our purposes.

Let $\Omega$ be a domain in $\bbR^n$ and let $U$ be an open set in $\bbR^n\times \bbR\times\bbR^n$
of the form $U = \Omega \times I \times B$, where $I$ is an open interval and $B$ is an open ball
in $\bbR^n$.  Consider the second order quasi-linear elliptic operator $Q=Q(u)$ as defined in Equation
\ref{eq:2d} for $U$-admissible functions $u\in C^2(\Omega)$, where  now we assume $a^{ij},b\in
C^{\8}(U)$.  We now also assume $Q=Q(u)$ is uniformly elliptic, by which we mean (1) the quantity 
$\sum_{i,j=1}^n a^{ij}(x,r,p)\xi^i\xi^j$ is uniformly positive and bounded away from
infinity for all $(x,r,p)\in U$ and all unit vectors $\xi=(\xi_1,...,\xi_n)$, and (2) $a^{ij},b$ and
their first order partial derivatives are bounded on $U$.     
  
We need the notion of a support function.  Given $u\in C^0(\Omega)$ and $x_0\in \Omega$, $\phi$
is an upper (resp., lower) support function for $u$ at $x_0$ provided $\phi(x_0) = u(x_0)$ and
$\phi \ge u$ (resp., $\phi \le u$) on some neighborhood of $x_0$.  We say that a function $u\in
C^0(\Omega)$ satisfies $Q(u)\ge 0$ \emph{in the sense of support
functions} iff for all $\e >0$ and
all $x\in \Omega$ there is a  $U$-admissible lower support function $\phi_{x,\e}$, which is 
$C^2$ in a neighborhood of $x$, such that
$Q(\phi_{x,\e})(x) \ge -\e$.  We say that $u$ satisfies $Q(u) \ge 0$ \emph{in the sense of support
functions with Hessians locally bounded from below} iff, in addition, there
is a constant $k>0$, independent of $\e$ and $x$, such that $\mathrm{Hess}(\phi_{x,\e})(x)\ge -kI$,
where $I$ is the identity matrix.     
For $u\in C^0(\Omega)$, we define $Q(u)\le 0$ \emph{in the sense of support functions} in an
analogous fashion.

The following theorem is a special case of Theorem 2.4 in \cite{AGH}

\begin{thm}\label{th:3e} 
Let $Q=Q(u)$ be a second order quasi-linear uniformly elliptic operator as described above.
Suppose  $u,v\in C^0(\Omega)$ satisfy, 
\begin{enumerate}
\item[(1)] $u\le v$ on $\Omega$ and $u(x_0)=v(x_0)$ for some $x_0\in \Omega$,  
\item[(2)] $Q(u) \ge 0$  in the sense of support functions with Hessians
locally bounded from below, and 
\item[(3)] $Q(v) \le 0$ in the sense of support functions.
\end{enumerate}
Then $u=v$ on $\Omega$ and $u = v \in C^{\8}(\Omega)$.
\end{thm}

We now proceed to the proof of the maximum principle for $C^0$ null hypersurfaces.

\smallskip\noindent
{\it Proof of Theorem \ref{th:3d}:} We first observe that $p$ is an interior point of 
null generators for both $S_1$ and $S_2$, and that these two null generators agree
near $p$.  To show this, let $\eta_i$ be a null generator of $S_i$ starting at $p$, $i=1,2$;
$\eta_1$ is future directed and $\eta_2$ is past directed.  Let $U$ be a neighborhood of $p$
in which $S_1$ is closed achronal , such that $S_2\cap U\subset J^+(S_1\cap U,U)$.  We may 
assume $\eta_1$ and $\eta_2$ are contained in $U$.  Since $\eta_2$ is past pointing,    
$\eta_2 \subset J^-(S_1\cap U,U)\cap J^+(S_1\cap U,U)$.  Then, as in Lemma \ref{th:3c}, the 
achronality of $S_1$ in $U$ forces $\eta_2\subset S_1$.  To avoid an achronality violation,
$\eta =-\eta_2 \cup \eta_1$ must be an unbroken null geodesic, and hence a null generator
of $S_1$ passing through $p$.  Similarly, $-\eta=-\eta_1\cup \eta_2$ is a null generator of
$S_2$ passing through $p$.   

Let $P$ be a timelike hypersurface passing through $p$ transverse to $\eta$.  Let $K_p$ be the
normalized semi-tangent of $S_1$ at $p$; $K_p$ is tangent to $\eta$.  Let $N$ be the spacelike
unit normal vector field to $P$ that points to the same side of $P$ as $K_p$.  By making a 
homothetic change in the background Riemannian metric we may assume $\la K_p,N_p\ra = 1$.  Hence,
$K_p$ is of the form, $K_p = Z_p+N_p$, where $Z_p\in T_pP$ is a future directed unit timelike
vector.  

As in the proof of Theorem \ref{th:2a}, $P$ in the induced metric can be expressed as in Equations
\ref{eq:2e} and \ref{eq:2f}.  Moreover, $V$ can be constructed so that $Z_p$ is perpendicular to $V$.
Then $K_p = (\d_0 + N)_p$.  By Lemma \ref{th:3b}, provided $P$ is taken small enough, $\S_1 = S_1\cap P$,
and $\S_2 = S_2\cap P$  will be partial Cauchy surfaces in $P$ passing through $p$, with $\S_2$ to the
future of $\S_1$.  Thus, shrinking $P$ further if necessary, there exist functions $u_i\in C^0(V)$, 
$i =1,2$, such that $\S_i = \mathrm{graph}(u_i)$ and

\begin{enumerate}
\item[(1)] $u_1\le u_2$ on $V$ and $u_1(p) = u_2(p) = 0$.
\end{enumerate}

Let $\{S_{q,\e}\}$ be the family of smooth null lower support hypersurfaces for $S_1$.
Restrict attention to those members of the family for which $q \in \S_1$.
Let $B_{q,\e}$ be the null second fundametal form of $S_{q,\e}$ at $q$ with respect to
the null vector $K_{q,\e}$.  By Lemmas \ref{th:3a} and \ref{th:3c}, the collection of null vectors
$\{K_{q,\e}\}$ can be made arbitrarily close to $K_p$ by taking $P$ sufficiently small. 
This has several implications. It implies, in particular, for $P$ sufficiently small, that
$K_{q,\e}$ is transverse to $P$. Hence, by shrinking $S_{q,\e}$  about $q$, if necessary,  
$S_{q,\e}$ meets $P$ in a properly transverse manner, and thus $\S_{q,\e}=S_{q,\e}\cap P$ is a
smooth spacelike hypersurface in $P$. Thus for each $\e>0$ and $q\in \S_1$, there exists $\phi_{q,\e}\in
C^{2}(W_{q,\e})$, $W_{q,\e}\subset V$, such that $\S_{q,\e} = \mathrm{graph}(\phi_{q,\e})$.

We now consider the null mean curvature operator $\th = \th(u)$, as described in equations
\ref{eq:2j} and \ref{eq:2jj}, on the set $U = V\times (-a,a)\times B$, where $B$ is an open
ball in $\bbR^{n-2}$ centered at the origin.  By choosing $V$, $a$ and $B$ sufficiently small,
$\th = \th(u)$ will be uniformly elliptic on $U$, in the sense described above.  
Since the vectors $\{K_{q,\e}\}$ can be made uniformly close to $K_p$,      
the inner products $\la K_{q,\e}, N_q\ra$ can be made uniformly close to the
value one.  Hence, we can rescale the vectors $K_{q,\e}$ so that $\la K_{q,\e}, N_q\ra=1$ without
altering the validity of the assumed null mean curvature inequality $\theta_1\ge 0$ in the sense of
support hypersurfaces,  at points of $S_1$ in $P$.  Then $K_{q,\e}$ can be expressed as,
$K_{q,\e} = Z_{q,\e}+N_q$, where $Z_{q,\e}\in T_qP$ is a future directed unit timelike vector.
Moreover, the vectors $Z_{q,\e}$ can be made uniformly close to $Z_p = \partial_0|p$ by taking  
$P$ small enough.  Equation \ref{eq:2hh} then implies that the Euclidean vectors 
$\d\phi_{q,\e}(q)=(\d_1\phi_{q,\e}(q),...,\d_{n-2}\phi_{q,\e}(q))$ can be made to lie in the ball $B$.

We conclude that by taking $P$ sufficiently small, each function $\phi_{q,\e}$
is $U$-admissible.  Now, $\phi_{q,\e}$ is a $C^2$ lower support function for $\S_1$
at $q$. By assumption, the null mean curvature of $S_{q,\e}$ at $q$ satisfies,
$\th_{q,\e}(q) \ge -\e$, which, in the  analytic setting, translates into,
$\th(\phi_{q,\e})(q)\ge -\e$.   Thus, $u_1$ satisfies,
$\th(u_1) \ge 0$ in the sense of support functions.

For each $q\in \S_1$, $Z_{q,\e}$ is the future directed timelike unit normal to
$\S_{q,\e}$.  Let $\beta_{q,\e}$ be the second fundamental form of $\S_{q,\e}\subset P$ at $q$
with respect to the normal $Z_{q,\e}$.  Let $B_{P,q}$ be the second fundamental form of P at $q$ with
respect to $N$.  The second fundamental forms $B_{q,\e}$, $\beta_{q,\e}$, and $B_{P,q}$ are related
by
\beq
B_{q,\e}(\ol X,\ol X) = \beta_{q,\e}(X,X) + B_{P,q}(X,X), \label{eq:3b}
\eeq
for all unit vectors $X\in T_q\S_{q,\e} = [Z_{q,\e}]^{\perp}\subset T_qP$. In a sufficiently small
relatively compact neighborhood $P_0$ of $p$ in $P$, the collection of vectors $\{Z_{q,\e}:q\in
\S_1\cap P_0\}$ has compact closure in $TP$.  It follows that the collection of vectors
$\mathcal{X}= \{X_q\in T_q\S_{q,\e}:  q\in \S_1\cap P_0, |X_q|=1\}$ has compact closure in $TP$, as
well. Hence the set of numbers $\{B_{P,q}(X_q,X_q): X_q\in \mathcal{X}\}$ is bounded.  Coupled with the
assumption that the second fundamenal forms $\{B_{q,\e}\}$ are locally bounded from below, we
conclude, by shrinking $P$ further if necessary, that the second fundamental forms
$\{\beta_{q,\e}: q\in \S_1\}$ are locally bounded from below, i.e., for each $q_0\in \S_1$ there is 
a neighborhood $\mathcal W$ of $q_0$ in $\S_1$ and a constant $k$ such that
\beq\label{eq:3c}
\beta_{q,\e} \ge- k g_{q,\e} \quad\mbox{for all } q\in \calW \mbox{ and } \e >0,
\eeq 
where $g_{q,\e}$ is the induced metric on $\S_{q,\e}$ at $q$.  
For $P$ sufficiently 
small, the induced metrics $g_{q,\e}$ will be close to the metric of $V$ at $p$.
Using the relationship between $\beta_{q,\e}$ and $\hess{\phi_{q,\e}}$, worked out, for 
example, in Section 3.1 in \cite{AGH},
inequality~\ref{eq:3c} translates into the analytic
statement that for each $q_0\in \S_1$ there is 
a neighborhood $\mathcal W$ of $q_0$ in $\S_1$ and a constant $k_1$ such that
$\hess{\phi_{q,\e}}(q) \ge -k_1I$ for all  $q\in \mathcal{W}$ and $\e >0$. 
Thus, we finally conclude that,
\begin{enumerate}
\item[(2)] $u_1$ satisfies
$\th(u_1) \ge 0$ in the sense of support functions  with Hessians locally bounded from below.
\end{enumerate}

\noindent
By similar reasoning, adjusting the size of $P$ as necessary, we have that
\begin{enumerate}
\item[(3)] $u_2$ satisfies
$\th(u_2) \le 0$ in the sense of support functions. 
\end{enumerate}
\noindent
In view of (1), (2), and (3), Theorem \ref{th:3e} applied to the operator $\th = \th(u)$ 
implies that
$u_1 =u_2$ on $V$ and $u_1 =u_2$ is $C^{\8}$.  Hence, $\S_1$ and $\S_2$ are smooth spacelike
hypersurfaces in $P$  which coincide near $p$.  Then near $p$, $S_1$ and $S_2$ are obtained by
exponentiating normally out along a common smooth null orthogonal vector field
along $\S_1=\S_2$.  The conclusion to 
Theorem \ref{th:3d} now follows.  
\hfill\endproof
\smallskip

For simplicity we have restricted attention to the null mean curvature inequalities
$\th_1\le 0 \le \th_2$.  However, Theorem \ref{th:3d}, with an obvious modification of 
Definition \ref{def:3b}, remains valid under the null mean curvature inequalities
$\th_1\le a \le \th_2$, for any $a \in \bbR$. 

\section{The null splitting theorem}

We now  consider some consequences of Theorem \ref{th:3d}.  The proof of
many global results in general relativity, such as the classical Hawking-Penrose
singularity theorems and more recent results such as those concerning topological censorship
(see e.g., \cite{FSW, GW})
involve the construction of a \emph{timelike line\/} or a \emph{null line\/}
in spacetime.  A timelike
geodesic segment is \emph{maximal\/} if it is
longest among all causal curves joining its end points, or equivalently, if it realizes the 
Lorentzian distance between its end points.
A timelike line is an inextendible timelike geodesic each segment of which is  maximal.
Similarly, a null line is an  inextendible null geodesic  each segment of which is
maximal.  But since each segment of a null geodesic has zero length, it follows that an
inextendible null geodesic is a null line iff it is achronal, i.e., iff no two points of it can
be joined by a timelike curve.  In particular, null lines, like timelike lines, must be free
of conjugate points.  The standard Lorentzian splitting theorem \cite[Chapter 14]{BEE},
considers the rigidity of spacetimes which contain timelike lines.  Recall, it asserts
that a timelike geodesically complete spacetime $(M,g)$ which obeys the strong energy
condition, $\mathrm{Ric\/}(X,X)\ge 0$ for all timelike vectors $X$, and which contains
a timelike line splits along the line, i.e., is isometric to 
$(\bbR\times V,-dt^2\oplus h)$, where $(V,h)$ is a complete Riemannian manifold.
We now consider the analogous problem for spacetimes with null lines.  The  theorem
stated below establishes the rigidity of spacetime in this null case. 
Unless otherwise
stated, we continue to assume that all null vectors are normalized to unit length with respect to
a fixed background Riemannian metric.  

\begin{thm}\label{th:4a} 
Let $M$ be a null geodesically complete spacetime which obeys the null energy condition,
$\mathrm{Ric}(X,X) \ge 0$ for all null vectors $X$ and contains a null line $\eta$.  Then $\eta$
is contained in a smooth closed achronal totally geodesic ($B\equiv 0$) null hypersurface $S$.  
\end{thm} 

The simplest illustration of Theorem~\ref{th:4a} is Minkowski space:  Each null geodesic in
Minkowski space is contained in a null hyperplane.  De Sitter space, anti-de Sitter, and gravitational
plane wave solutions furnish other illustrations.

\begin{rem}\label{rem:4a}{\rm
Theorem \ref{th:4a} may be viewed as a ``splitting" theorem of sorts, where
the splitting takes place in $S$.
Let $K$ be the unique (up to scale) future
pointing null vector field of $S$.
The vanishing of the null second fundamental form $B$
of $S$  means that the standard kinematical quantities associated with $K$,
i.e., the expansion $\theta$, and shear $\sigma$ (as well as vorticity) 
of $K$ vanish.  In a similar vein,  
the vanishing of $B$ implies that the metric $h$ on $TS/K$ defined in 
Section 2.1 is invariant under the flow generated by $K$; see \cite{K} for
a precise statement and proof. To take this a step further, suppose, in the setting
of Theorem \ref{th:4a}, that $M$ is also globally hyperbolic.  Let $V$ be a smooth Cauchy surface
for $M$.  Since $V$ is Cauchy, and the generators of $S$ are inextendible in $S$, each null generator of
$S$ meets $V$ exactly once, and this intersection is properly transverse.  Hence,
$\S=S\cap V$ is a smooth codimension two spacelike submanifold of $M$, and the map $\Phi:\bbR\times \S\to
S$, defined by,
$\Phi(s,x) = \exp_xsK(x)$, is a diffeomorphism.  The invariance of the metric $h$ with
respect to the flow generated by $K$, or equivalently, by $\Phi_*(\frac{\d}{\d s})$, implies that
$(i\circ \Phi)^*g = \pi^*g_0$,
where $i: S \hookrightarrow M$ is inclusion, $\pi : \bbR\times \S \to \S$ is projection, and $g_0$
is the induced metric on $\S$.  In  more heuristic terms, $S\approx \bbR \times \S$,
and $i^*g = 0dt^2 + g_0$.  It is in this sense that one may view Theorem~\ref{th:4a} as a splitting
theorem.} 
\end{rem}  

\begin{rem}\label{rem:4b}
{\rm  The proof of Theorem \ref{th:4a} shows that $S = \widehat{\d I^+(\eta)} = \widehat{\d I^-(\eta)}$, 
where $\widehat{\d I^{\pm}(\eta)}$ is the component of $\d I^{\pm}(\eta)$ containing $\eta$.
In more heuristic terms $S$ is obtained as a limit of future null cones $\d I^+(\eta(t))$
(resp., past null cones $\d I^-(\eta(t))$) as $t\to -\8$ (resp., $t\to \8$).  The proof
also shows that the assumption of null completeness can be weakened. It is sufficient to
require that
the generators of $\d I^-(\eta)$ be future geodesically complete and the generators of 
$\d I^+(\eta)$ be past geodesically complete. (As discussed in the proof, the generators of 
$\d I^-(\eta)$ (resp., $\d I^+(\eta)$) are in general future (resp., past)
inextendible in $M$.) }
\end{rem}

The proof of Theorem~\ref{th:4a} is an
application of  Theorem \ref{th:3d}, and relies on the following basic lemma.

\begin{lem}\label{th:4b} Let $M$ be a spacetime which satisfies the null energy condition.  Suppose $S$
is  an achronal $C^0$ future null hypersurface whose null generators are future geodesically complete.
Then $S$ has null mean curvature $\theta\ge 0$ in the sense of support hypersurfaces, with null second 
fundamental forms locally bounded from below.
\end{lem} 

\begin{proof} 
The basic idea for constructing past support hypersurfaces for $S$ is to consider the ``past
null cones" of points on generators of $S$ formed by past null geodesics emanating from
these points.  For the purpose of establishing certain properties about these null hypersurfaces,
it is useful to relate them to achronal boundaries of the form $\d J^-(q)$, $q\in S$, which
are defined purely in terms of the causal structure of $M$.  For this reason
we assume initially that $M$ is globally hyperbolic.  At the end we will indicate
how to remove this assumption.  

For each $p \in S$ and normalized semi-tangent $K$ at $p$, let
$\eta : [0,\8)\to M$, $\eta(s) = \exp_psK$, be the affinely parameterized null geodesic generator
of $S$ starting at $p$ with initial tangent $K$.  Since $S$ is achronal each such generator
is a null \emph{ray}, i.e., a maximal null half-geodesic.  
Since $\eta$ is maximal, no point on $\eta|_{[0,r)}$
is conjugate to $\eta(r)$.
For each $r>0$, consider the 
achronal boundary $\d J^-(\eta(r))$, which is a $C^0$ hypersurface containing $\eta|_{[0,r]}$.
By standard properties of the null cut locus (see especially, Theorems 9.15 and 9.35 in \cite{BEE}, which
assume global hyperbolicity) there is a neighborhood $U$ of $\eta|_{[0,r)}$ such that  $S_{p,K,r}=\d
J^-(\eta(r))\cap U$ is a smooth  null hypersurface diffeomorphic under the exponential map based at
$\eta(r)$ to a neighborhood of the line segment $\{-\tau\eta'(r): 0< \tau\le 1\}$ in the
past null cone $\Lambda_{\eta(r)}^-\subset T_{\eta(r)}M$. 

From the achronality of $S$ we observe, $\d J^-(\eta(r)) \cap I^+(S)= \emptyset$.
This implies that $S_{p,K,r}$ is a past support hypersurface for $S$ at $p$.  
Let $\th_{p,K,r}$ denote the null mean curvature of  $S_{p,K,r}$ at $p$ 
with respect to $K$.
We use Equation \ref{eq:2c} to obtain the lower bound
\beq\label{eq:4a}
\th_{p,K,r} \ge - \frac{n-2}{r}.
\eeq
The argument is standard.  In the notation of Section 2.1, let $\th(s)$, $s\in [0,r)$,
be the null mean curvature of $S_{p,K,r}$ at $\eta(s)$ with respect to $\eta'(s)$.
 Equation \ref{eq:2c} and the energy condition imply,
\beq\label{eq:4b}
\frac{d\th}{ds} \le - \frac1{n-2}\,\th^2.
\eeq 
Without loss of generality we may assume $\th(0) = \th_{p,K,r} <0$.  Then $\th = \th(s)$ is
strictly negative on $[0,r)$, and we can devide \ref{eq:4b} by $\th^2$ to obtain,
\beq\label{eq:4c}
\frac{d}{ds}\, \th^{-1} \ge \frac1{n-2}.
\eeq
Integrating \ref{eq:4c} from $0$ to $r-\delta$, and letting $\delta\to 0$ we obtain the lower
bound~\ref{eq:4a}. 

Thus, since $r$ can be taken arbitrarily large, we have shown, in the globally hyperbolic case, that $S$
has null mean curvature $\th\ge 0$ in the sense of support hypersurfaces with respect to the collection
$\{S_{p,K,r}\}$ of smooth null hypersurfaces.  

By Lemma \ref{th:3c}, $K$ is tangent to $S_{p,K,r}$ at $p$.
Let $B_{p,K,r}$ denote the null second fundamental form of $\{S_{p,K,r}\}$ at $p$ with respect
to $K$.  We now argue that the collection of null second fundamental forms $\{B_{p,K,r}: r\ge r_0\}$,
for some $r_0 >0$,
is locally bounded from below.  Recall, ``locally", means ``locally in the
point $p$"; the lower bound must hold for all $r > r_0$, cf., inequality \ref{eq:3bb}.
This lower bound follows from a continuity argument and an elementary monotonicity
result, as we now discuss.

Fix $p\in S$.  Let $U$ be a convex normal neighborhood of $p$. Thus, for each $q\in U$,
$U$ is the diffeomorphic image under the exponential map based at $q$ of a neighborhood
of the origin in $T_qM$.  Provided $U$ is small enough, we have that for each $q\in U$,
$\d J^-(q) \cap U = \d(J^-(q) \cap U) = \exp_q(\Lambda_q^-\cap \exp_q^{-1}(U))$, where $\Lambda_q^-$ is the
past null cone in $T_qM$.  In particular, for each $q\in U$, $W(q) = \d J^-(q) \cap U \setminus \{q\}$
is a smooth null hypersurface in $U$ such that if $q_n\to q$ in $U$, $W(q_n)$ converges smoothly
to $W(q)$.  

There exists a neighborhood $V$ of $p$, $V\subset U$, and $r_0 >0$ such that for each $q\in V$,
and each normalized null vector $K\in T_qM$, the null geodesic segment $s\to \eta(s) = \exp_qsK$,
$s\in [0,r_0]$ is contained in $U$.  Let $B(q,K)$ be the null second fundamental form of $W(\exp_qr_0 K)$
at $q$ with respect to $K$.  Since, as $(q_n,K_n) \to (q,K)$, $W(\exp_{q_n}r_0 K_n)$ converges smoothly
to $W(\exp_qr_0 K)$, we have that  $B(q_n,K_n)\to B(q,K)$ smoothly.  Returning to the original
family of support hypersurfaces, $\{S_{p,K,r}\}$, with associated family of null
second fundamental forms $\{B_{p,K,r}\}$, observe that when $q\in S\cap V$, 
$S_{q,K,r_0}$ agrees with $W(\exp_qr_0 K)$ near $q$.  Hence, if $(q_n,K_n) \to (q,K)$ in $S\cap V$,
$B_{q_n,K_n,r_0}\to B_{q,K,r_0}$ smoothly.  It follows that the collection of null second fundamental
forms $\{B_{p,K,r_0}\}$ is locally bounded from below.  

Consider as in the beginning, for $p\in S$ and $K$ a normalized semi-tangent at $p$,
the null geodesic generator, $\eta : [0,\8)\to M$, $\eta(s) =
\exp_psK$, $s\ge 0$.  For $0 < r < t < \8$, $J^-(\eta(r)) \subset J^-(\eta(t))$. Then,
since $\d J^-(\eta(t))$ is achronal, $\d J^-(\eta(r))$ cannot enter $I^+(\d J^-(\eta(t))$.
It follows that for $r<t$, $S_{p,K,r}$ lies to the past side of $S_{p,K,t}$ near $p$.
By  an elementary comparison of null second fundamental forms at $p$ we obtain
the monotonicity property, 
\beq\label{eq:4cc}
B_{p,K,t}\ge B_{p,K,r} \quad \mbox{ for all }0 < r < t < \8.
\eeq
This monotonicity now implies that the entire family
of null second fundamental forms $\{ B_{p,K,r}: r \ge r_0\}$ is locally bounded from
below.  

This concludes the proof of Lemma \ref{th:4b} under the assumption that $M$ is globally hyperbolic.
We now describe how to handle the general case.  
In general, $M$ may have bad causal properties. In particular the generators of $S$ could be closed.
Nevertheless, the past support hypersurfaces are formed in the same manner, as 
the ``past
null cones" of points on generators of $S$ formed by past null geodesics emanating from
these points.  But as an intermediary step, to take advantage of the arguments in the
globally hyperbolic case, we pull back each generator 
to a spacetime having good causal properties.

Again, consider, for $p\in S$ and $K$ a normalized semi-tangent at $p$,
the null geodesic generator, $\eta : [0,\8)\to M$, $\eta(s) =
\exp_psK$, $s>0$.  Restrict attention to the finite segment $\eta|_{[0,r]}$.
Roughly speaking, we introduce Fermi type coordinates near $\eta|_{[0,r]}$.
Let  $\{e_1, e_2, ... e_{n-1}\}$ be an orthonormal frame of spacelike vectors
in $T_{\eta(0)}M$.  Parallel translate these vectors along $\eta$ to obtain
spacelike orthonormal vector fields $e_i = e_i(s)$, $1 \le i \le n-2$ along
$\eta|_{[0,r]}$.    Consider the map $\Phi : \bar M \subset \bbR^n\to M$ defined by 
$\Phi(s,x^1, x^2, ..., x^{n-1}) =\exp_{\eta(s)}(\sum_{i=1}^{n-1}x^ie_i)$.  Here $\bar M$
is an open set containing the line segment $\{(s,0,0, ...,0): 0\le s\le r\}$.  By
the inverse function theorem  we can choose $\bar M$ so that $\Phi$ is a local
diffeomorphism.  Equip $\bar M$ with the pullback metric $\bar g= \Phi^*g$,
where $g$ is the Lorentzian metric on $M$, thereby making $\bar M$ Lorentzian and $\Phi$ a local
isometry.  Let $t\in C^{\8}(\bar M)$ be the $0$th
coordinate function, $t(s,x^1, x^2, ..., x^{n-1})=s$.  Since the 
slices $t=s$ are spacelike, $\D t$ is
timelike, and hence $t$ is a time function on $\bar M$.  Thus, $\bar M$ is a strongly causal
spacetime.  

The curve $\bar \eta:[0,\bar r+\delta]\to \bar M$, $\bar\eta(s) =(s, 0, 0, ... , 0)$,
defined for $\delta>0$ sufficiently small,
is a maximal null geodesic in $\bar M$.  Let $\mathcal{K}$ be a compact neighborhood of $\bar \eta$ in
$\bar M$. Then by Corollary 7.7 in \cite{P2}, any two causally related points in $\mathcal{K}$ can be
joined by a longest causal curve $\gamma$ in $\mathcal{K}$, and each segment of $\gamma$ contained in the
interior of $\mathcal{K}$ is a maximal causal geodesic.  This property is sufficient to push through, with
only minor modifications, all relevant results concerning the null cut locus of $\bar\eta(0)$ on $\bar
\eta$.    In view of the above discussion,
there is a neighborhood $\bar U$ of $\bar\eta|_{[0,r)}$ such that  $\bar S_{p,K,r}=\d
J^-(\eta(r))\cap\bar U$ is a smooth  null hypersurface
diffeomorphic under the exponential map based at
$\bar\eta(r)$ to a neighborhood of the line segment $\{-\tau\bar\eta'(r): 0< \tau\le 1\}$ 
in the past null cone $\Lambda_{\bar\eta(r)}^-\subset T_{\bar\eta(r)}\bar M$.
Let $\bar V\subset \bar U$ be 
a neighborhood of $\bar
\eta(0)$ on which $\Phi$ is an isometry onto its image, and let   
$S_{p,K,r}=\Phi(\bar S_{p,K,r}\cap\bar V)$.
Then $\{S_{p,K,r}\}$  is the desired collection of past support
hypersurfaces for
$S$, having all the  requisite properties.  In particular, the mean curvature inequality \ref{eq:4a} and
the monotonicity property \ref{eq:4cc} hold for the family $\{\bar S_{p,K,r}\}$, and hence the family
$\{S_{p,K,r}\}$, by just the same arguments as in the globally hyperbolic case.  
\end{proof}

\smallskip\noindent
{\it Proof of Theorem \ref{th:4a}:}  Since $\eta$ is achronal, it follows that
$\eta \subset \partial I^-(\eta)$, and hence $\partial I^-(\eta) \ne \emptyset$.
Then, by standard properties of achronal boundaries, $\d I^-(\eta)$
is a closed achronal $C^0$ hypersurface in $M$.  We claim that $\d I^-(\eta)$ is a $C^0$ future
null hypersurface.  
By standard results on achronal boundaries, e.g., \cite[Theorem~3.20]{P2},
$\partial
I^-(\eta)\setminus \bar \eta$ (where $\bar \eta$ is the closure of $\eta$ as a subset
of $M$) is ruled by future inextendible null geodesics.  However, since we 
do not assume $M$ is strongly causal, it is possible, in the worst case scenario, that
$\bar\eta = \partial I^-(\eta)$, in which case \cite[Theorem 3.20]{P2} gives no information.  To show that
$\partial I^-(\eta)$ is ruled by future inextendible null geodesics we apply instead the 
more general \cite[Lemma~3.19]{P2}.
Let $N\subset U$ be a 
convex normal neighborhood of $p$.  $N$ as a spacetime in its own right is strongly causal.
Let $K$ be a compact neighborhood of $p$ contained in $N$.  Then for each $t\in \bbR$, $\eta|_{[t,\8)}$  
cannot remain in $K$ if it ever meets it.  Thus there exists a sequence of $p_i = \eta(t_i)$,
$t_i\nearrow \8$, $p_i\notin K$.  It follows that for each $x\in K\cap I^-(\eta)$, there exists a
future directed timelike curve from $x$ to a point on $\eta$ not in $K$.
We may then apply \cite[Lemma 3.19]{P2} to conclude that $p$ is the past end point of a 
future directed null geodesic segment contained in $\partial I^-(\eta)$.
Hence, according to Definition \ref{def:3a},
$\d I^-(\eta)$ is a $C^0$ future null hypersurface. 
Moreover, because $\d I^-(\eta)$ is closed,
the  null generators of $\d I^-(\eta)$ are future inextendible in $M$, and hence future complete.  

Let $S_-$ be the component of
$\d I^-(\eta)$ containing $\eta$.  From the above, $S_-$ is 
an achronal $C^0$ future null hypersurface whose null generators are future 
geodesically complete.  Thus, by Lemma \ref{th:4b}, $S_-$ has null mean curvature 
$\theta_-\ge 0$ in the sense of support hypersurfaces, with null second 
fundamental forms locally bounded from below.  Let $S_+$ be the component 
of $\d I^+(\eta)$ containing $\eta$.  Arguing in a time-dual fashion,
$S_+$ is 
an achronal $C^0$ past null hypersurface whose null generators are past 
geodesically complete.  By the time-dual of Lemma \ref{th:4b}, 
$S_+$ has null mean curvature 
$\theta_+\le 0$ in the sense of support hypersurfaces.  Moreover at any point
$p$ on $\eta$, $S_+$ lies to the future side of $S_-$ near $p$.  Theorem~\ref{th:3d}
then implies that $S_- =S_+ =S$ is a smooth  null hypersurface
containing $\eta$ with vanishing null mean curvature, $\theta \equiv 0$.   
Equation \ref{eq:2c} and the null energy condition
then imply that the shear scalar vanishes along each generator, which in turn implies that the null second
fundamental form of $S$ vanishes. 
\hfill\endproof
\smallskip

We conclude the paper with an application of Theorem \ref{th:4a}.   The  application we consider is
concerned with asymptotically simple  (e.g., asymptotically flat and nonsingular) spacetimes as 
defined by Penrose \cite{P1} in terms of the notion of conformal infinity.

Consider a $4$-dimensional connected chronological spacetime $M$ with metric $g$ which can be conformally
included into a spacetime-with-boundary $M'$ with metric $g'$ such that $M$ is the
interior of $M'$, $M = M'\setminus \d M'$.  Regarding the conformal factor,
it is assumed that there exists a smooth function $\Omega$ on $M'$ 
which satisfies, 
\begin{enumerate}
\item[(i)] $\Omega > 0$ and $g'=\Omega^2g$ on $M$, and
\item[(ii)] $\Omega = 0$ and $d \Omega \ne 0$ along $\d M'$.
\end{enumerate}

The boundary $\scri := \d M'$ is assumed to consist of two components, $\scri^+$ and $\scri^-$,
future and past null infinity, respectively, which are smooth null hypersurfaces.  
$\scri^+$ (respectively,
$\scri^-$) consists  of points of $\scri$ which are future (resp., past) endpoints of causal curves in $M$.
A spacetime $M$ satisfying the above is said to be \emph{asymptotically flat at null infinity}.
If, in addition, $M$ satisfies,
\begin{enumerate}
\item[(iii)] Every inextendible null geodesic in $M$ has a past end point on $\scri^-$ and
a future end point on $\scri^+$
\end{enumerate}
then $M$ is said to \emph{asymptotically simple} with null conformal boundary.  
Condition (iii) is imposed
to ensure that $\scri$ includes all of the null infinity of $M$. 
It also implies  that $M$ is free of singularities which would
prevent some null geodesics from reaching infinity. 

The notion of asymptotic simplicity was introduced by Penrose  
in order to facilitate the study of the asymptotic behavior of
isolated gravitating systems.  When restricting to   
vacuum (i.e., Ricci flat) spacetimes, asymptotic simplicity
provides an elegant and rigorous setting for the study of 
gravitational radiation far from the radiating source. 
See the papers \cite{F2, F3} of
Friedrich for  discussion of the problem of global existence 
of asymptotically simple vacuum spacetimes.
Here we prove the following rigidity result.  

\begin{thm}\label{th:4c} 
Suppose $M$ is an asymptotically simple vacuum spacetime which
contains a null line.  Then $M$ is isometric to Minkowski space.  
\end{thm} 

\begin{proof} The first and main step of the proof is to show that $M$ is flat (i.e.,
has vanishing Riemann curvature).  Then certain global arguments will show
that $M$ is isometric to Minkowski space.

For technical reasons, it is useful to extend $M'$ slightly beyond its boundary. 
In fact,
by smoothly attaching a collar to $\scri^+$ and to $\scri^-$, we can extend
$M'$ to a spacetime $P$ without boundary such that $M'$ is a retract of $P$ and both 
$\scri^+$ and $\scri^-$ separate $P$.  It follows that $\scri^+$ and $\scri^-$ are globally
achronal null hypersurfaces in $P$.  

Let $M^- = M\cup \scri^-$.  A straight forward limit curve argument, using the asymptotic simplicity
of $M'$,
shows that $M^-$ is causally simple.  This means that the sets of the form $J^{\pm}(x,M^-)$,
$x\in M^-$, are closed subsets of $M^-$. 
(The limit curve lemma, in the form of Lemma 14.2
in \cite{BEE}, for example,
is valid in~$P$.)  

Let $\eta$ be a null line in $M$ which has past end point $p\in\scri^-$ and
future end point $q\in \scri^+$.  
Consider the ``future null cone" at $p$, $N_p :=\d I^+(p,M^-)$.  From the causal
simplicity of $M^-$ it follows that $\d I^+(p,M^-) = J^+(p,M^-) \setminus I^+(p,M^-)$.
Hence each point in $N_p$  can be joined to $p$
by a null geodesic segment in $M^-$.  In particular $N_p$ is connected. From the
simple equality $I^+(p, M^-) =I^+(\eta, M^-)$, it follows that, 
$$N_p = \d I^+(\eta, M^-) = \d I^+(\eta, M) \cup \g_p = S\cup \g_p\, ,$$ 
where
$S =\d I^+(\eta, M)$ and $\g_p$ is the future directed null generator of $\scri^-$ starting
at $p$.

Since asymptotically simple spacetimes are null geodesically complete, Theorem \ref{th:4a}
implies that $S= \d I^+(\eta, M)$ is a smooth null hypersurface in $M$ which is
totally geodesic with respect to $g$.  Since it is smooth and closed the generators of $S$ never
cross and never leave $S$ in $M^-$ to the future.  Moreover, since  $N_p=\d I^+(p,M^-)$
 is achronal, there are no conjugate  points to $p$ along
the generators of $N_p$.  It follows that $N_p\setminus \{p\}$ is the diffeomorphic
image under the exponential map $\exp_p: T_pP\to P$ of $(\Lambda_p^+\setminus\{0\})\cap \exp_p^{-1}(M^-)$, 
where $\Lambda_p^+$ is the future null cone  in $T_pP$.  We are now fully
justified in referring to
$N_p$ as the future null cone in $M^-$ at $p$. 

Given a smooth null hypersurface, with smooth future pointing null normal
 vector field $K$, the shear tensor $\sigma_{ab}$ is the trace free part of the
null second fundamental form $B_{ab}$, $\sigma_{ab}=B_{ab} -\frac{\theta}{2}h_{ab}$.
Since $S$ is totally geodesic in the physical metric $g$ and the shear 
scalar $\s=(\s_{ab}\s^{ab})^\frac12$ is a conformal 
invariant, it follows by continuity that the shear tensor of $N_p\setminus \{p\}$
(with respect to an appropriately chosen $g'$- null normal $K'$) vanishes
in the metric $g'$, $\sigma_{ab}'= 0$.  The trace free part of equation \ref{eq:2b}
then implies that the components $C_{a0b0}'$
(with respect to an appropriately chosen pseudo-orthonormal frame
in which  $e_0=K'$, cf., Section 4.2 in \cite{HE}) 
of the conformal tensor of $g'$ vanish on $N_p\setminus\{p\}$.  
An argument of Friedrich in \cite{F1}, which makes use of the
Bianchi identities and, in the present case, the vacuum field equations expressed in
terms of his regular conformal field equations,  then shows that
the so-called rescaled conformal tensor,
and hence, the conformal tensor of the 
physical metric $g$ must vanish on $D^+(N_p,M^-)\cap M$.  
Since $M$ is Ricci flat, we conclude that $M$ is
flat on $D^+(N_p,M^-)\cap M$.  In a precisely time-dual fashion $M$ is 
flat on $D^-(N_q,M^+)\cap M$, where
$N_q$ is the past null cone of $M^+= M\cup\scri^+$ at $q$.

To conclude that $M$ is everywhere flat we show that $M\subset D^+(N_p,M^-)\cup D^-(N_q,M^+)$.
Consider the set $V=I^+(S,M)\cup S\cup I^-(S,M)$.  It is clear from
the fact that $S$ is an achronal boundary  that $V$ is
open in $M$. As $S= \d I^+(\eta,M)$, it follows that $I^+(S,M) \subset J^+(p,M^-)$, and  
time-dually, that $I^-(S,M)\subset J^-(q,M^+)$.  
Using the fact that $J^+(p,M^-)$ and $J^-(q,M^+)$ are closed subsets 
of $M^-$ and $M^+$, respectively, it follows that $V$ is closed in $M$.  Hence, $M=I^+(S,M)\cup S\cup
I^-(S,M)$. We show that each term in this union is a subset of $D^+(N_p,M^-)\cup D^-(N_q,M^+)$.  Trivially,
$S\subset N_p \subset D^+(N_p,M^-)$.  
Consider $I^+(S,M) \subset J^+(p,P)\cap M=J^+(N_p,P)\cap M$. 
We claim that  $J^+(N_p,P)\cap M\subset D^+(N_p,P)\cap M$.  If not, then  
$H^+(N_p,P)\cap M\ne\emptyset$. Choose a point $x\in H^+(N_p,P)\cap M$, and let $\nu$ be a null
generator of $H^+(N_p,P)$ with future end point $x$.  
Since $N_p$ is edgeless, $\nu$ remains in $H^+(N_p,P)$
as it is extended into the past.  By asymptotic simplicity, $\nu$ must meet
$\scri^-$.  In fact, since  no portion of $\nu$ can coincide with a generator of $\scri^-$, $\nu$ must
meet $\scri^-$ transversely and then enter $I^-(\scri^-,P)$.  But this means that $\nu$ has
left $J^+(N_p,P)$, which is a contradiction.  Since  $D^+(N_p,P)\cap M =D^+(N_p,M^-)\cap M$,
we have shown that $I^+(S,M)\subset D^+(N_p,M^-)$.  By the time-dual argument, $I^-(S,M)\subset
D^-(N_q,M^+)$.

Thus, $M$ is globally flat.  Also, as an asymptotically simple spacetime,
$M$ is  null geodesically complete, simply connected and globally hyperbolic; cf., 
\cite{HE}, \cite{N}.  We use these properties to show that $M$ is geodesically
complete.  It then follows from the  uniqueness of 
simply connected space forms that $M$ is isometric
to Minkowski space. 

We first observe that  there exists a time orientation preserving
local isometry $\phi:M \to L$, where $L$ is Minkowski space.  To obtain $\phi$, first construct
by a standard procedure a frame $\{e_0,e_1,...,e_{n-1}\}$ of orthonormal parallel vector fields
on $M$.  Then solve $dx^i= \la e_i,\,\ra$, $i= 0,...,n-1$,
for functions $x^i\in C^{\8}(M)$, and set $\phi = (x^0,x^1,...,x^{n-1})$. 

From the fact that $\phi$ is a local isometry
and $M$ is null geodesically complete, it follows that any null geodesic, or broken
null geodesic, in $L$ can be lifted via $\phi$ to $M$.  In particular it follows that
$\phi$ is onto:  If $\phi(M)$ is not all of $L$ then we can find a null geodesic $\bar\eta$
in $L$ that meets $\phi(M)$ but is not entirely contained in $\phi(M)$.  Choose $p\in M$
such that $\phi(p)$ is on $\bar\eta$.  Then the lift of $\bar\eta$ through $p$ is incomplete,
contradicting the null geodesic completeness of $M$.      

We now show that $M$ is timelike geodesically complete.  If it is not, then, without loss
of generality, there exists a future inextendible unit speed timelike geodesic
$\g:[0,a)\to M$, $t \to \g(t)$, with $a < \8$.  Let $\bar\g = \phi\circ\g$; $\bar\g$ can be extended
to a complete timelike geodesic in $L$ which we still refer to as $\bar\g$.    
Let $\bar\eta$ be a future directed broken null geodesic extending from $\bar p=\bar\g(0)$ 
to $\bar q=\bar\g(a)$.  Let $\eta$ be the lift of $\bar\eta$ starting at $p = \g(0)$; $\eta$
extends to a point $q\in I^+(p)$, with $\phi(q) = \bar q$.  Since $M$ is globally hyperbolic
there exists a maximal timelike geodesic segment $\mu$ from $p$ to $q$.  Then $\bar\mu =\phi\circ\mu$
is a timelike geodesic segment in $L$ from $\bar p$ to $\bar q$.  Hence, $\bar\mu = \bar\g|_{[0,a]}$.
It follows that $\mu$ extends $\g$ to $q$, which is a contradiction. 

Finally, we show that $M$ is spacelike geodesically complete.  
If it is not, then there exists an
inextendible unit speed spacelike geodesic
$\g:[0,a)\to M$, $t \to \g(t)$, with $a < \8$. 
Let $\bar\g = \phi\circ\g$; $\bar\g$ can be extended
to a complete spacelike geodesic in $L$ which we still refer to as $\bar\g$.  We now use the
fact that timelike geodesics, and broken timelike geodesics, in $L$ can be lifted via $\phi$
to $M$.  Let $\bar\a=\bar\a_1+\bar\a_2$ be a two-segment broken timelike geodesic extending 
from $\bar p=\bar\g(0)$ to $\bar q=\bar\g(a)$, with $\bar\a_1$ starting at $\bar p$ and past pointing.
Similarly, let $\bar\b=\bar\b_1+\bar\b_2$ be a two-segment broken timelike geodesic extending 
from $\bar p$ to $\bar q$, with $\bar\b_1$ starting at $\bar p$ and future pointing.
Let $\bar x$ be the point at the corner of $\bar \a$, and let $\bar y$ be the point at the corner 
of $\bar \b$.  Note that $\bar\g|_{[0,a]}\subset I^+(\bar x) \cap I^-(\bar y)$.
Let $\a$ be the lift of $\bar \a$ starting at $p= \g(0)$; $\a$ extends to a point $q$ with 
$\phi(q) = \bar q$.  Similarly, let $\b$ be the lift of $\bar \b$ starting at $p$. 
Let $x$ be the point at the corner of $\a$, and let $y$ be the point at the corner 
of $\b$.  

Note that an initial segment of $\g$ is contained in  $J^+(x)\cap J^-(y)$.
We claim that $\g$ is entirely contained in $J^+(x)\cap J^-(y)$.  If not then
$\g$  either leaves $J^+(x)$ or $J^-(y)$.  Suppose it leaves $J^+(x)$ at $z\in \g \cap\d J^+(x)$.  Since
$M$ is globally hyperbolic, there exists a null geodesic segment $\eta$ from $x$ to $z$. Then
$\bar \eta = \phi\circ \eta$ is a null geodesic in $L$ from $\bar x$ to $\bar z \in \bar\g|_{[0,a]}$.
But since, by construction, $\bar\g|_{[0,a]}\subset I^+(\bar x)$, no such null geodesic exists.

Thus, $\g\subset J^+(x)\cap J^-(y)$.  We show that $\g$ extends
to $q$, thereby obtaining a contradiction. Consider a sequence $\g(t_n)$, $t_n\to a$.  Since $J^+(x)\cap
J^-(y)$ is compact, there exists a subsequence $\g(t_{n_k})$ which converges to a point $q'\in J^+(x)\cap
J^-(y)$.  Since $\bar\g(t_n)\to \bar q$, it follows by continuity that $\phi(q')=\bar q$.  Let
$\mu$ be a causal geodesic segment from $x$ to $q'$.  Then $\bar \mu =\phi\circ\mu$ is a causal
geodesic from $\bar x$ to $\bar q$ in $L$.  Thus, $\bar \mu = \bar\a_2$, and hence $\mu=\a_2$.
Since $\a_2$ has future end point $q$, we conclude that $q'=q$.  Hence, since every sequence $\g(t_n)$,
$t_n\to a$, has  a subsequence converging to $q$, it follows that $\lim_{t\to a}\g(t)= q$, and so $\g$
extends to
$q$. This concludes the proof that $M$ is geodesically complete and hence, as noted above, 
isometric to Minkowski space. 
\end{proof} 

We remark in closing that Theorem \ref{th:4c} can be generalized in various directions.  
For example,  it is possible to formulate a version of Theorem \ref{th:4c} for asymptotically
flat spacetimes which contain singularities, black holes, etc., by imposing suitable conditions
on the domain of outer communications $D= I^-(\scri^+)\cap I^+(\scri^-)$, the conclusion then
being that $D$ is flat.  Also, it appears that Theorem~\ref{th:4c} can be extended to vacuum 
spacetimes with positive cosmological constant, $\Lambda > 0$, thereby yielding a characterization
of de Sitter space.  Details of this will appear elsewhere. 

\smallskip
\noindent
\textbf{Acknowledgements.}
The author is indebted to Helmut Friedrich for many helpful discussions concerning the proof
of Theorem \ref{th:4c}.  The author would also like to thank Lars Andersson and
Piotr Chru\'sciel for many valuable comments and suggestions.  Part of the work on this
paper was carried out during a visit to the Royal Institute of Technology in Stockholm, Sweden
in 1999. The author wishes to thank the Institute for its hospitality and support.

\providecommand{\bysame}{\leavevmode\hbox to3em{\hrulefill}\thinspace}

\end{document}